\documentclass[12pt]{article}
\usepackage{epsfig}
\usepackage{epsf}
\usepackage{here}
\usepackage{color}
\usepackage{pifont}

\newcommand{\dpg}[2]{\frac{\partial #1}{\partial #2}}
\newcommand{\ddpg}[2]{\frac{\partial^2 #1}{\partial #2^2}}

\newcommand{\Ob}{{{\cal O}}}

\newcommand{\I}{{{\cal I}}}
\newcommand{\LL}{{\cal L}}

\newcommand{\Et}{{\rm E}_{\theta}}
\newcommand{\Etet}{{\rm E}_{\theta_*}}
\newcommand{\var}{{\rm var}}

\begin{document}

\title{  A NEWTON-LIKE ALGORITHM FOR LIKELIHOOD MAXIMIZATION: THE ROBUST-VARIANCE SCORING ALGORITHM} \author{{\sc Daniel Commenges, H\'el\`ene Jacqmin-Gadda,} \\ {\sc C\'ecile Proust and J\'er\'emie Guedj}\\ INSERM E03 38, ISPED\\ Universit\'e Victor Segalen Bordeaux 2\\146 rue L\'eo
Saignat, Bordeaux, 33076, France\\Tel: (33) 5 57 57 11 82; Fax (33) 5 56 24 00 81}

\maketitle

\vspace{\baselineskip} {\bf }
This article considerss a Newton-like method already used by several authors but which has not been thouroughly studied yet. We call it the robust-variance scoring (RVS) algorithm because the main version of the algorithm that we consider replaces minus the Hessian of the loglikelihood used in the Newton-Raphson algorithm by a matrix $G$ which is an estimate of the variance of the score under the true probability, which uses only the individual scores. Thus an iteration of this algorithm requires much less computations than an iteration of the Newton-Raphson algorithm. Moreover this estimate of the variance of the score estimates the information matrix at maximum. We have also studied a convergence criterion which has the nice interpretation of estimating the ratio of the approximation error over the statistical error; thus it can be used for stopping the iterative process whatever the problem. A simulation study confirms that the RVS algorithm is faster than the Marquardt algorithm (a robust version of the Newton-Raphson algorithm); this happens because the number of iterations needed by the RVS algorithm is barely larger than that needed by the Marquardt algorithm while the computation time for each iteration is much shorter. Also the coverage rates using the matrix $G$ are satisfactory.

\vspace{3ex} Keywords: likelihood, Newton algorithm, convergence tests, score, robustness. \vspace{10mm}

\section{INTRODUCTION}
The Newton-Raphson method is very efficient in maximizing the likelihood in problems of moderate complexity such as the Cox model for instance. 
 However, the development of more and more complex statistical models in particular random-effect models and  semi-parametric models together with more complex observation schemes is a challenge for likelihood inference. In such complex situations we may not be in a convex optimization problem and we need more robust methods than the Newton-Raphson methods; moreover the likelihood is difficult to compute, most often involving non-analytical integrals,  and analytical derivatives are not available. A robust version of the Newton-Raphson method can be developed using the idea of  the Marquardt (or Levenberg-Marquardt) algorithm (Marquardt, 1963), but this algorithm needs the second derivatives of the loglikelihood (the Hessian). 

These difficulties have prompted the development of alternative algorithms. The most often used for likelihood inference is the EM algorithm (Dempster, Laird and Rubin, 1977). Many researchers also have adopted a Bayesian framework in which the MCMC algorithm can be used (Gilks et al., 1996). There are also applications of the MCMC algorithm for likelihood maximization using the EM algorithm (Wei and Tanner, 1990). However both EM and the Bayesian approaches may be very time-consuming in complex problems, although it is an active research field to improve the algorithms using these approaches: for instance Kuhn and Lavielle (2005) have proposed a fast version of the EM algorithm. Newton-like methods still have trumps in their hands. The first thing to note is that numerical derivatives can be easily computed; there are also semi-analytical ways of computing the derivatives in cases where analytical derivatives are available for a ``full'' problem and this has been used for instance by Hedeker and Gibbons (1994) and Commenges and Rondeau (2006). Advantages of Newton-like methods are the relatively small number of iterations for converging, the existence of good convergence tests and the inference that can be done using an estimator of the information matrix at maximum. However the computation of the Hessian needed in the Newton-Raphson and Marquardt algorithm, even if feasible, may be time-consuming and inaccurate.

The aim of this paper is to study a Newton-like algorithm proposed by Berndt et al. (1974); the idea can be traced back to Bock (1970) and the algorithm can be viewed as an extension of the Gauss-Newton algorithm. It has been extensively used in the psychometric and econometrics literature where it is called the BHHH algorithm, and maybe less in the biostatistical literature for which we may cite  Hedeker and Gibbons (1994), Todem, Kim and Lesaffre (2006) and Ruppert (2005). We feel that it  has still not been studied thoroughly: in particular far from the maximum the definition of the scoring matrix should be modified (Meilijsson, 1989). This algorithm estimates the information matrix at maximum by an estimator of the variance of the score using individual scores and far from the maximum the matrix is still an estimate of the variance of the score; for this reason we call it the Robust-variance scoring (RVS) algorithm. This algorithm requires much less computations than the Marquardt algorithm while retaining all of its advantages. Moreover we argue that a very good convergence test is available for this algorithm and that robust inference can be naturally done in case of mis-specified models.

The organization of the paper is as follows. In section 2 we present some facts that are valid for a class of Newton-like algorithms; in particular we propose a convergence test which has an intrinsic statistical interpretation. In section 3 we present the RVS algorithm in detail and we show that robust inference can be done without computing the Hessian. 
In section 4 we make some numerical experiments using a random-effect model with a complex observation scheme. Section 5 concludes.

\section{NEWTON-LIKE METHODS FOR LIKELIHOOD MAXIMIZATION}
\subsection{NOTATIONS}
We assume a model $\{P_{\theta}\}_{\theta \in \Theta}$, where $\Theta$ is a ``nice'' open subset of $\Re^m$ so that this defines a regular parametric model (Rothenberg, 1971; Bickel et al, 1993). Let $\LL_{\Ob}^{\theta}$ be the likelihood for the observation depicted by a $\sigma$-field $\Ob$; we shall assume that the observation is made of $n$ independent parts so that $\Ob=\vee \Ob_i$ where the $\Ob_i$ are independent $\sigma$-fields, for instance the $\Ob_i$ are generated by independent random variables $Y_i$. We note $L^{\theta}_i=\log \LL^{\theta}_{\Ob_i}$, the log-likelihood for observation $i$; $L^{\theta}=\sum_{i=1}^n L^{\theta}_i$ the total log-likelihood. We shall define the maximum likelihood estimator as the solution of minimizing $-L^{\theta}$; this has the advantage of removing some minus signs in many equations; moreover most optimization problems are presented as minimization problems and the MLE can be viewed as minimizing a Kulback-Leibler divergence. We assume that the log-likelihood is twice-differentiable and we note $U^{\theta}_i=-\dpg {L^{\theta}_i}{\theta}$, $U^{\theta}=\sum_{i=1}^n U^{\theta}_i$ and $H(\theta)=-\ddpg{L^{\theta}}{\theta}$. The information matrix is $I(\theta)=\Et [H(\theta)]= \Et [U^{\theta}(U^{\theta})^T]$. We assume that the $L^{\theta}_i$'s are independently identically distributed (iid) (the iid case is more general than may appear at first sight): then we have $I(\theta)=n\I(\theta)$, where $\I(\theta)$ does not depend on $n$. We assume that  $\I(\theta)$ is positive definite for all $\theta \in \Theta$. For the statistician there is a true probability measure $P_*$; if $P_*$ belongs to the model there is a value $\theta_*\in \Theta$ such that $P_*=P_{\theta_*}$. We shall assume first that this is the case (we shall consider later that this is not the case, that is that the model is mis-specified). Under these conditions $\theta_*$ is identifiable and the maximum likelihood estimator $\hat \theta$ is consistent and asymptotically efficient and we have $n^{-1/2}(\hat \theta -\theta_*) \rightarrow N(0,\I(\theta_*)^{-1})$.
\subsection{A CLASS OF NEWTON-LIKE METHODS}
We shall consider a class of Newton-like methods for likelihood maximization. The Newton algorithm is defined by the iteration 
$$\theta_{k+1}=\theta_k-H^{-1}(\theta_k)U(\theta_k).$$
This is the most efficient method for minimizing or maximizing a function when it is not too far from quadratic. When the function to maximize is the log-likelihood an additional advantage of the method is that $H(\hat \theta)$ is an estimator of $I(\theta_*)$ allowing to construct confidence intervals and tests for $\theta$. In complex problems however this method has two main drawbacks: it may fail, in particular because the Hessian may not be positive-definite; it may be time consuming if the number of parameters is large and the likelihood difficult to compute. So we are led to consider Newton-like methods defined by iterations of the type:
$$\theta_{k+1}=\theta_k-\alpha_kG^{-1}(\theta_k)U(\theta_k),$$
where $\alpha_k$ is the step length and is found by line-search and 
where $G(\theta_k)$ has two main  properties:

i) it is positive-definite for all $\theta$;

ii) it approaches $H(\theta_k)$ and $I(\theta_k)$ when $\theta_k$ approaches $\hat \theta$.

Another property is important:

iii) if the model is re-parametrized in $\eta_k=A\theta_k+a$ the following relation holds between $G$ and the matrix $G_{\eta}$ for the transformed problem: $G(\theta_k)=A^TG_{\eta}(\eta_k)A$.

Moreover the ease of computation of  $G(\theta_k)$ will be a major criterion for choosing among the methods which are successful in converging.

Property (i) ensures that the algorithm is a descent method and thus will converge;  property (ii) ensures that convergence near the maximum will be approximately as fast as with the Newton-Raphson algorithm and  that inference about $\theta$ will be possible by estimating $I(\theta_*)$ by $G(\hat \theta)$. 
Property (iii) ensures that the algorithm is invariant by affine transformation of the parameters, that is we should not have problems with bad conditioned $G$, as long as $G$ is non-singular and the numerical computations sufficiently accurate (see Fletcher, 1987). The three properties will also help finding a good convergence test.

\subsection{THREE NEWTON-LIKE METHODS}
Let us consider three Newton-like algorithms: the Marquardt, the Fisher-scoring and the RVS algorithms. 
The Marquardt algorithm (Marquardt, 1963) uses a matrix $G$ obtained by adding to $H$ a positive-definite matrix. The simplest form is $G(\theta_k)=H(\theta_k)+\lambda_k I_d$ where $I_d$ is the identity matrix and $\lambda_k$ is adaptively tuned so as to make the matrix $G(\theta_k)$ positive-definite, ensuring to $G(\theta_k)$ the property (i); if $H$ is positive definite around the maximum, smaller and smaller $\lambda$ are used so that $G$ for Marquardt has property (ii). The Marquardt algorithm has been found to be efficient in many statistical problems (for instance Joly et al., 2002; Alioum et al., 2005; Proust and Jacqmin-Gadda, 2005), although this is a general method. The computational load for computing $G(\theta_k)$ is about the same as for computing $H(\theta_k)$.

The Fisher-scoring algorithm uses $G(\theta_k)=I(\theta_k)$. With our assumptions, we have property (i) and obviously property (ii). This choice also enjoys property (iii). The Fisher scoring algorithm has been used in particular in generalized linear models where $I(\theta_k)$ is generally easy to compute; it is an efficient and robust algorithm (McCullagh and Nelder, 1999). The algorithm is specific to likelihood maximization. In some particular cases $I(\theta_k)$ is easy to compute but this is not the case in general. In general problems, when there is no analytical form of the information matrix, the Fisher-scoring algorithm seems difficult to implement.

Berndt et al. (1974) (see also Hedeker and Gibbons, 1994) used $G(\theta_k)=\sum_{i=1}^n U_i(\theta_k) U_i^T(\theta_k)$. Todem, Kim and Lesaffre (2006) used $G(\theta_k)=\sum_{i=1}^n U_i(\theta_k) U_i^T(\theta_k)-n^{-1}U(\theta_k) U^T(\theta_k)$, a correction proposed by Meilijsson (1989).  However none of these authors have thoroughly studied this type of algorithm. Note that the two choices are equivalent near the maximum since $U(\hat \theta)=0$ but not far from the maximum. Both choices enjoy the three properties (i), (ii) and (iii), although it may happen that (ii) is not satisfied with the second choice. The algorithm converged satisfactorily in the statistical problems treated by these authors. The main advantage of this type of choice is that $G(\theta_k)$ is much less computationally demanding than $H(\theta_k)$ since it requires only computing the $U_i(\theta_k)$'s, that is only the first derivatives of the individual log-likelihoods. This advantage becomes huge when the number of parameters is large. This is the type of algorithm that we will study in section 4 under the name of ``Robust-variance scoring'' (RVS) algorithm.

\section{CONVERGENCE TESTS FOR NEWTON-LIKE METHODS}
\subsection{CONVENTIONAL CONVERGENCE TESTS}

Using a good  convergence test is essential in the efficiency of any iterative method. We shall call ``convergence test'' a method leading to the decision of stopping: generally it will take the form $C_k<c$, where $C_k$ will be called ``stopping criterion'' and $c$ stopping value.
Stopping criteria based on the displacement in the parameter space $||\theta_{k+1}-\theta_k||^2$ or in the log-likelihood $L^{\theta_{k+1}}-L^{\theta_k}$ are not satisfactory because small displacements may occur if the algorithm fails to find a good direction. A better criterion is based on $||U(\theta_k)||^2$ because a necessary condition for the maximum is $U(\hat \theta)=0$. It is nice to have a stopping criterion which is invariant by linear transform of $\theta$; such a stopping criterion is based on $U(\theta_k)G^{-1}(\theta_k)U(\theta_k)$, where $G$ has property (iii) (Fletcher, 1987: Dennis and Schnabel, 1996). It still remains to fix a value at which to stop. 

We can say a little more on the criterion $U(\theta_k)G^{-1}(\theta_k)U(\theta_k)$ in the context of log-likelihood maximization.
Dennis and Schnabel (1996) have noted that this stopping criterion is not invariant by a change of scale of the function to maximize; however the log-likelihood does not have to be rescaled: it is defined up to an additive constant (which corresponds to the choice of the reference probability when defining the likelihood) but it can not be multiplied by an arbitrary constant without changing its meaning.
Another property, also specific to likelihood maximization, is that the criterion is asymptotically invariant near the maximum by any one-to-one transforms. This follows from Delta-calcul.

\subsection{AN OPTIMAL CONVERGENCE TEST}

Here we shall develop an ``optimal'' stopping criterion for maximizing the likelihood.
 It is good to have a stopping criterion which has the invariance properties mentioned above. However it is still difficult to choose a value at which to stop. A way to solve the problem is to devise a criterion which can be interpreted and which has the same meaning for different problems. The proposed solution is to stop when the approximation error is small with regard to the statistical error. The approximation error for a given $\hat \theta$ can be measured by $d(\theta_k,\hat \theta)$, where $d(.,.)$ is a distance defined on $\Re^m$. The statistical error can be defined as $\Etet [d(\hat \theta, \theta_*)]$, so that the criterion is
$$C_k=\frac{d(\theta_k,\hat \theta)}{\Etet [d(\hat \theta, \theta_*)]}$$
and it is natural to take the stopping value $c$ between $10^{-2}$ and $10^{-4}$, whatever the problem. 
Let us define the distance $d(x,y)$ by a norm of $x-y$ of the type: $d(x,y)=(x-y)^TM(x-y)$. It is desirable that the distance be invariant by an affine transformation of the parameters, so it is natural to take $M=G$, the matrix used in our algorithm and which has property (iii). Near the maximum we have by property (ii) that $G(\theta_k) \approx I(\theta_*)\approx I(\hat \theta)$. Thus we have $\Etet [d(\hat \theta, \theta_*)]\approx \Etet [(\hat \theta- \theta_*)^TI(\theta_*)(\hat \theta- \theta_*)]$. Using the asymptotic property of $\hat \theta$ the right-hand term is the expectation of a chi-squared with $m$ degrees of freedom; thus  $\Etet [d(\hat \theta, \theta_*)]\approx m$. As for the numerator, using the first order Taylor expansion $U(\theta_k)\approx I(\hat \theta)(\theta_k-\hat \theta)$ it can be seen that we have $d(\theta_k,\hat \theta) \approx U(\theta_k)G^{-1}(\theta_k)U(\theta_k)$. Thus we have an approximation of the criterion $C_k$ as:
$$C_k \approx m^{-1}U(\theta_k)^TG^{-1}(\theta_k)U(\theta_k).$$
This is a scaled version of the invariant criterion mentioned in the previous paragraph, but now we have an interpretation of it and we can choose a stopping value which does not need to depend on the problem at hand.

A good stopping criterion must not take small values far from the maximum, a criticism that has been addressed to other criteria in section 3.1. It is difficult to prove that $C_k$ has this property for general $G$ because the only properties that we have specified and which are relevant far from the maximum, (i) and (iii), are not sufficient to impose large values of $C$. We shall see in section 4.3 that the result is easy to obtain for the RVS algorithm. Also in the appendix we put forward the idea that in certain difficult problems we might decide to stop the iterations before full convergence.

\section{THE ROBUST-VARIANCE SCORING (RVS) ALGORITHM}
\subsection{DESCRIPTION OF THE BASIC ALGORITHM}
Let us study the algorithm based on $G(\theta_k)=\sum_{i=1}^n U_i(\theta_k) U_i^T(\theta_k)-\eta_k n^{-1}U(\theta_k) U^T(\theta_k)$. Hedeker and Gibbons (1994) chose $\eta_k=0$. This choice may lead to inefficient directions; it can be seen that it leads to small displacements for the components having larges values of $U$. Todem, Kim and Lesaffre (2006) chose $\eta_k=1$. The above problem does not appear with this choice. Moreover $G(\theta_k)$ can then be interpreted as an estimator of the variances of $U$ under the true model, and this warrants the name of ``Robust-variance scoring algorithm''. Indeed, under certain assumptions, the law of large number gives us: $n^{-1}G(\theta_k) \rightarrow n^{-1}\var _{P_*}[U(\theta_k)]$.  Note that $\var _{P_*}[U(\theta_k)]\neq \var _{\theta_k}[U(\theta_k)]=I(\theta_k)$, so that the RVS algorithm is different from the Fisher-scoring algorithm, even asymptotically.

One possible problem with the choice $\eta_k=1$ is that it is not sure that $G(\theta_k)$ is not singular. One possibility to ensure property (i) is to use $\eta_k$ slightly smaller than $1$, possibly in an adaptive way as in the Marquardt algorithm; an adaptive choice of $\eta_k$ can also be done for optimizing the direction of search.

\subsection{INFERENCE USING $G$}
We may use $G(\hat \theta)$ as an estimator of $\var _{P_*}[U(\theta_*)]$. This may be slighlty biased when the number of parameters is large relatively to $n$  and one might consider reducing the bias by using $\frac{n}{n-m}G(\hat \theta)$. If the model is well specified we have $\var _{P_*}[U(\theta_*)]=\var _{\theta_*}[U(\theta_*)]=I(\theta_*)$, so that tests or confidence intervals based on $G(\hat \theta)$ are asymptotically equivalent to tests based on $H(\hat \theta)$. So we may estimate $\var _{\theta_*}(\hat \theta)$ by $G^{-1}(\hat \theta)$.

 If the model is mis-specified (that is $P_*$ does not belong to $\{P_{\theta}\}_{\theta \in \Theta}$), $\theta_*$ has the meaning that $P_{\theta_*}$ is the probability measure of $\{P_{\theta}\}_{\theta \in \Theta}$ which is the closest to $P_*$ with respect to Kullback-Leibler divergence (Burnham and Anderson, 1998) and $\hat \theta$ is consistent for $\theta_*$. Royall and Tsou (2003) call $\theta_*$ the object of inference: for inference to be relevant with a mis-specified model the object of inference must have the same meaning in the true probability $P_*$. When the model is mis-specified we no longer have $\var _{P_*}[U(\theta_*)]=I(\theta_*)$; inference based on $I(\theta_*)$ is not valid while robust inference based on $G(\hat \theta)$ in the spirit of Royall (1986) is still feasible. The variance of the estimators of the parameters can be obtained by the sandwich estimator:
$\widehat \var _{P_*}(\hat \theta)= H^{-1}(\hat \theta)G(\hat \theta)H^{-1}(\hat \theta)$. However this estimator has the disadvantage of requiring the computation of the Hessian while the algorithm itself avoided it. In some cases this is feasible to compute it once for the computation of confidence intervals. 

In the case where it is not feasible to compute the Hessian it should be possible to construct a robust confidence interval using the following idea. For any value $\theta$ a score test of ``$\theta_*=\theta$'' can be done: the null hypothesis is rejected if $U^T(\theta)G^{-1}(\theta)U(\theta)>c_{\alpha}$, where we use the $\chi ^2_m$ asymptotic distribution of the statistic and $c_{\alpha}$ is the critical value for constructing a size-$\alpha$ test based on this distribution. Thus the set $\{\theta : U^T(\theta)G^{-1}(\theta)U(\theta)\le c_{\alpha}\}$ is a confidence ellipsoid at level $\alpha$. It seems possible to derive algorithms for constructing confidence interval along these lines.

\subsection{CONVERGENCE TEST}
The convergence test described in section 3.2 can be applied to the RVS algorithm and we have seen that it enjoys a very nice interpretation near the maximum. We wish in addition that the probability of stopping before reaching the region of the maximum be very low. 

Note that with our assumption that $I(\theta)>0$ there is a unique $\theta$ such that $\Etet [U(\theta)]=0$ and this is $\theta_*$; in case where the model is well specified this value defines the true probability measure. This does not mean that there is a unique maximum of the likelihood for given observation but the probability of several maxima tends toward zero. Suppose that we wish to test the hypothesis $\Etet [U(\theta_k)]=0$; a score test statistic is $U(\theta_k)^T[\var _{P_*}[U(\theta_k)]]^{-1}U(\theta_k)$; since $G(\theta_k)$ is an estimator of $\var _{P_*}[U(\theta_k)]$, $C_k$ can be considered as a test statistic for $\Etet [U(\theta_k)]=0$. If $c$ is the stopping value there is identity between the decision of not stopping and the decision of rejecting the null hypothesis. Thus our convergence test is a statistical test (outside of the region of the maximum). The null hypothesis is rejected if $U(\theta_k)^TG(\theta_k)^{-1}U(\theta_k)>cm$; as $c$ is small the size of the test will be large. This is of course not a conventional testing situation and what we wish is that the power be large. If $\Etet [U(\theta_k)]\neq 0$ the distribution of $U(\theta_k)^TG(\theta_k)^{-1}U(\theta_k)$ is asymptotically a non-central chi-squared and the probability of rejecting the null hypothesis tends toward $1$ (that is, the probability of stopping is close to zero for large $n$). In fact even if we neglect the non-centrality term, the probability of stopping is very low with the choice $c=10^{-2}$ and $m$ moderately large or large. For instance for $m=10$ we have $P_*(C_k<c)<P(\chi^2_{10}<0.1)=2.5 10^{-8}$.

\subsection{SPEED}
The algorithm is useful for problems in which the computation of the likelihood is difficult and time-consuming and the derivatives of the loglikelihood can not be computed analytically. In that case the derivatives will be computed by finite difference (Overton, 1981) and the computational load for one iteration is approximately proportional to the number of likelihood evaluations for computing $U(\theta_k)$ and $G(\theta_k)$. For computing $U(\theta_k)$ one has to compute the $U_i(\theta_k)$, so that $G(\theta_k)$ is a by-product of the computation of $U(\theta_k)$ for the RVS algorithm. This computation takes $2m$ evaluations of the loglikelihood if one uses centered differences for approximating the derivative. 
This is to be compared to the Newton or Marquardt algorithm for which one has to compute $H(\theta_k)$ which requires the computation of $m(m+1)/2$ second derivatives. Even if these second derivatives are computed with only one likelihood evaluation this leads to $m(m+5)/2$ evaluations. Thus the ratio of computation required for one iteration of the RVS compared to the Marquardt algorithm is approximately $4/(m+5)$ which takes values $0.26, 0.16, 0.11, 0.08$ for $m= 10, 20, 30, 40$ respectively; that is for large problems the computational burden for one iteration of the RVS algorithm is about ten times less than for the Marquardt algorithm. 

In terms of accuracy, iterating the numerical differentiation, necessary for the Marquardt algorithm, may lead to unacceptable loss of precision by cancellation error (Overton, 1981) although  in some cases there is the possibility to make more precise computations (Commenges and Rondeau, 2006). For the RVS algorithm the numerical differentiation is not iterated, diminishing the accuracy problem; in some cases it is also possible to avoid the numerical differentiation. A problem may occur in the RVS algorithm if $n$ is small because the rationale of the algorithm is based on the fact that $G(\theta_k)$ is an estimator of $\var _{P_*}[U(\theta_k)]$ using the law of large numbers; as mentioned by Song et al. (2005) this estimator may be unstable for small $n$ jeopardizing both convergence of the algorithm and inference based on $G(\hat \theta)$.

\subsection{APPLICATION TO PENALIZED LIKELIHOOD}
Penalized likelihood is defined as:
\begin{equation} \label{penlik} pl(\theta)= L^{\theta} -J(\theta,\kappa),\end{equation}
where $\kappa=( \kappa_j, j=1,\ldots, K)$ is a set of smoothing coefficients.
Penalized likelihood is useful for obtaining smooth non-parametric estimators of functions (O'Sullivan, 1986; Joly et al., 2002), so that $\theta$ may contain functions. However the solution is generally approximated on a basis of splines so that we are driven back to the parametric case, where part of the parameters are splines coefficients. 

For applying the RVS algorithm it is tempting to use the same formula as above replacing the $U_i$'s by $\dpg{pl_i}{\theta}$ where $pl_i=L^{\theta}-n^{-1}J(\theta,\kappa)$. However for the choice $\eta_k=1$ we then obtain the same matrix as for the ordinary likelihood because this matrix estimates $\var (U)$ and $ \var (\dpg{pl}{\theta})=\var (U)$. If we wish to have an efficient direction near the maximum it is appropriate to use:
$$G(\theta_k)=\sum_{i=1}^n U_i(\theta_k) U_i^T(\theta_k)-\eta_k n^{-1}U(\theta_k) U^T(\theta_k)-\ddpg{J(\theta,\kappa)}{\theta}.$$ 
With this choice $G(\theta_k)\approx \ddpg{pl}{\theta}(\theta_k)$ near the maximum. It is generally easy to compute $\ddpg{J(\theta,\kappa)}{\theta}$.

\section{NUMERICAL EXPERIMENTATION}
We performed a simulation study to compare RVS with the Marquardt algorithm for computing maximum likelihood estimators of linear mixed models from left censored longitudinal data. Indeed left censoring of biological measures (such as HIV-RNA)  frequently arises due to the lower detection limit of the measurement tools. Jacqmin-Gadda et al. (2000) have shown that the contribution to the likelihood for one subject when some measures are left censored is the product of the multivariate Gaussian density for the completely observed measures and of a multivariate Gaussian distribution function for the censored measures given the observed ones. Computation of the distribution function for each subject requires a numerical integration of size equals to the number of censored measures. Thus, the computation time highly depends on the proportion of censored measures and thus on the threshold value.

Data were generated according to either a linear model with an auto-regressive error structure or a linear model with random intercept and slope:

Model AR :
$$Y_{ij}= \beta_0 +  \beta_1 t_{ij} + \beta_2 X_{i} +  \beta_3 X_{i}t_{ij} + w_{ij} + e_{ij}$$

with $\beta_0=4$,  $\beta_1=-0.5$, $\beta_2=-0.5$, $\beta_3=-0.1$. The Gaussian variable $w_{ij}$ has mean 0, and covariance $cov(w_{ij},w_{il})=\sigma^2_w exp(- \delta |t_{ij}-t_{il}|)$ with $\sigma^2_w=1$,  $\delta=0.1$. The variance of the independent Gaussian random error $\sigma^2_e$ was 1. 

Model RE :
$$Y_{ij}= \beta_0 +  \beta_1 t_{ij} + \beta_2 X_{i} +  \beta_3 X_{i}t_{ij} + a_{0i}+ a_{1i}t_{ij} + e_{ij}$$

The two random effects $a_{0i}$ and $a_{1i}$ have a centered Gaussian distribution with 
 $\sigma^2_0=var(a_{0i})=0.25$, $\sigma^2_1=var(a_{1i})=0.1$ and $\sigma_{01}=cov(a_{0i},a_{1i})=-0.1$. The rest of the model was unchanged.

The response values less than the threshold were left censored. The threshold was either $1.0$ or $2.0$ leading respectively to about 25\% and 45\% of censored measures.

For each model and each threshold, we generated 100 samples of 100 subjects with a number of measures randomly selected according to a uniform distribution between 5 and 11. The covariate X had a Bernoulli distribution with P=0.5 and the time variable followed a U[0-6] distribution.

The starting values for the estimation algorithms were the maximum likelihood estimates obtained by imputing the threshold for the censored measures. Convergence was reached when the stopping criteria $C_k$ defined in section 3.2 reached the stopping value 10$^{-4}$ in less than 30 iterations for the Marquardt algorithm and in less than 30*(m+5)/4 for RVS (with m=7 or 8 according to the model). Results were not different when increasing the maximum number of iterations. In both algorithms, derivatives were computed by finite difference and a line-search was used to select the step length $\alpha_k$ only when the updated parameters did not improve the likelihood because this strategy was found to be efficient.

Convergence results are displayed in Table 1. The proportion of successful convergence was always higher for the RVS algorithm and the convergence time was clearly lower even if the number of iterations required to reach convergence was slightly higher. The ratio of computation time required for one iteration of RVS compared to the Marquardt algorithm was very close to the expected value 4/(m+5) and even slightly more favorable to RVS than expected (0.29 vs 0.33 expected for the AR model) because the Marquardt algorithm requires more computations in addition to computations of the derivatives (diagonal inflation and maybe more often line search for the step length).

When both algorithms converged the estimates were nearly identical: the mean sum of squared differences between the  estimates was less than 3.10$^{-6}$ in the four situations. Table 2 presents the mean asymptotic variance, the sample variance and the coverage rate of the 95\% confidence interval of the estimates for the two models and the two algorithms when the censoring threshold was $2.0$. This shows that the variances of the estimators were correctly estimated by the two methods ($H^{-1}$ or $G^{-1}$). Only two coverage rates were significantly different from the nominal values 95\% (with the AR model and the Marquardt algorithm).

\section{CONCLUSION}
We have studied a Newton-like method for maximizing the likelihood which needs only the computation of the scores and does not need the computation of the Hessian. 
A simulation study of random effects models with left censored observations confirms that the RVS algorithm is faster than the Marquardt algorithm and this happens because the number of iterations needed by the RVS algorithm is not much larger. The advantage of the RVS algorithm over methods needing computation of the Hessian increases with the number of parameters. In fact the Levenberg-Marquardt idea of increasing the diagonal of the matrix for improving its condition number could be applied to the RVS algorithm. The simulation study showed that the coverage rate using the matrix $G$ was good. We have also studied a convergence criterion for stopping the iterative process which has a statistical interpretation leading to an easy choice of the stopping value. 
Further work is needed for studying inference in misspecified models without computing the Hessian and this could be developed along the lines of section 4.2.
The RVS is the algorithm of choice for maximum likelihood inference in complex problems when the Hessian of the loglikelihood is difficult to compute and when the number of parameters is large. It can be applied in particular for estimating parameters in models specified by systems of non-linear differential equations (Guedj et al., 2007). 
\section*{REFERENCES}
\noindent
\setlength{\parindent}{-8mm}

Bickel, P; J., Klaassen, C. A. J., Ritov, Y. and Wellner, J. A. (1993).  {\it Efficient and Adaptive Estimation for Semiparametric Models}. New-York: Springer-Verlag.\vspace{3mm}

Berndt, E.,  Hall, B., Hall, R., Hausman, J. (1974). Estimation and Inference in Nonlinear Statistical Models. {\em Annals of Economic and Social Measurement} {\bf 3}, 653-665.\vspace{3mm}

Bock, R.D. and Lieberman, M. (1970). Fitting a response model forn dichotomously scored items. {\em Psychometrika} {\bf 35}, 179-197.\vspace{3mm}

Burnham, K. P., and D. R. Anderson.  (1998).  {\em Model selection and inference:
a practical information-theoretic approach}. Springer-Verlag, New York.\vspace{3mm}

Commenges, D. and Rondeau, V. (2006). Relationship between derivatives of the observed and full loglikelihoods and application to algorithms.The International Journal of Biostatistics: Vol. 2: No. 1, Article 4. http://www.bepress.com/ijb/vol2/iss1/4.\vspace{3mm}

Dennis, J. E. and Schnabel R. B. (1996). {\em Numerical methods for unconstrained optimization and nonlinear
equations}. Prentice-Hall (Englewood Cliffs, NJ).\vspace{3mm}

Dempster, A.D., Laird, N.M. and Rubin, D.B. (1977). Maximum likelihood from incomplete data via the EM
algorithm. {\em J. R. Statist. Soc. B} {\bf 39}, 1.\vspace{3mm}

Fletcher, R. (1987). {\em Practical Methods of Optimization}. John Wiley \& Sons (Chichester).\vspace{3mm}

Gilks, W. R., Richardson, S. and Spiegelhalter, D. J. (1996). {\em Markov Chain Monte Carlo in practice.}
New-York: Chapman \& Hall.\vspace{3mm}

Guedj, J., Thi\'ebaut, R. and Commenges, D. (2007). Maximum likelihood estimation in dynamical models of HIV. {\em Biometrics}, in press.\vspace{3mm}

Hedeker D. and Gibbons R.D. (1994). A random-effects ordinal regression model for multilevel analysis. {\em
Biometrics} {\bf 50}, 933-944.\vspace{3mm}

Jacqmin-Gadda, H., Thi\'ebaut, R., Ch\^ene, G., Commenges, D. (2000).
Analysis of left-censored longitudinal data with application to viral load in HIV infection. {\em Biostatistics} {\bf 1}, 355-368.\vspace{3mm}

Joly, P., Letenneur, L., Alioum, A., and Commenges, D. (1999). PHMPL: a computer program for hazard estimation using a penalized likelihood method with interval-censored and left-truncated data. {\em Computer Methods and Programs in Biomedicine} {\bf 60}, 225-231\vspace{3mm}

Kuhn, E. and Lavielle, M. (2005). Maximum likelihood estimation in ninlinear mixed effects models. {\em Computational Statistics \& Data Analysis} {\bf 49}, 1020-1038.\vspace{3mm}

Knight, K. (2000). {\em Mathematical Statistics}. Chapman \& Hall.\vspace{3mm}

Marquardt, D. (1963). An algorithm for least-squares estimation of nonlinear parameters. {\it SIAM Journal of
Applied Mathematics} {\bf 11}, 431-441.\vspace{3mm}

Meilijson, I. (1989). A fast improvement to the $\mathrm{EM}$ algorithm on its own terms. {\em Journal of the Royal Statistical Society. Series B } {\bf 51}, 127-138.\vspace{3mm}

Overton, M. (2001). {\em Numerical computing with IEEE floating point arithmetic}.  Philadelphia:
SIAM.\vspace{3mm}

Proust, C. and Jacqmin-Gadda, H. (2005).
 Estimation of linear mixed models with a mixture of distribution for the random-effects. {\em Computer Methods and Programs in Biomedicine} {\bf 78},
  165-173.\vspace{3mm}

Royall, R.M. (1986)  Model robust confidence intervals using maximum likelihood estimators
{\em International Statistical Review} {\bf 54 }, 221-226. \vspace{3mm}

Royall, R.M.  and Tsou, T.S. (2003).(Interpreting statistical evidence by using imperfect models: robust adjusted likelihood functions.{\em Journal of the Royal Statistical Society: Series B}
{\bf 65}, 391-404.  \vspace{3mm}

Ruppert, D.. (2005). Comment on ``Maximization by Parts in Likelihood Inference''
{\em Journal of the American Statistical Association} {\bf 100}, 1161-1163. \vspace{3mm}

Song, P.X., Fan, Y.  and  Kalbfleisch, J.D. (2005). Maximization by Parts in Likelihood Inference 
{\em Journal of the American Statistical Association} {\bf 100}, 1145-1158. \vspace{3mm}

Todem, D. , Kim, K-M and Lesaffre, E. (2006).
Latent-Variable Models for Longitudinal Data with Bivariate Ordinal Outcomes. {\em Statistics in Medicine}, in press.\vspace{3mm}

Wei, G.C.G. and Tanner, M.A. (1990).
Monte Carlo Implementation of the EM Algorithm and the Poor Man's Data Augmentation Algorithms
{\em Journal of the American Statistical Association} {\bf 85}, 699-704.\vspace{3mm}

\vspace{3mm}

\setlength{\parindent}{-0mm}
\section*{APPENDIX: STOPPING BEFORE FULL CONVERGENCE}
In some difficult problems it may save time to stop before full convergence, at least in an exploratory phase. We compute $\tilde U= U+ \varepsilon $ where $\varepsilon$ is a numerical error. The matrix $G$ is also affected by numerical error but there is a specific problem with $U$ because the norm of $U$ decreases near the maximum, so that the relative error increases. The computed search direction is $G^{-1}(\theta_k) \tilde U(\theta_k)=G^{-1}(\theta_k) U(\theta_k)+G^{-1}(\theta_k) \varepsilon_k$. The $G$-norm of the first term is $U(\theta_k)^TG^{-1}(\theta_k)U(\theta_k) \approx d(\theta_k,\hat \theta)$ and for the second term it is $\varepsilon ^TG^{-1}(\theta_k)\varepsilon$, so that the relative error in the search direction is 
$r_e=\frac{\varepsilon ^TG^{-1}(\theta_k)\varepsilon}{d(\theta_k,\hat \theta)}$ which increases as $\theta_k$ approaches $\hat \theta$. In particular when $r_e\approx 1$ we do not have one exact significant digit so that further progression toward $\hat \theta$ becomes impossible. This means that in term of stopping value it will be nearly impossible to go below $\tilde c=\frac{\varepsilon ^TG^{-1}(\theta_k)\varepsilon}{m}$. In difficult problems (where the computation of the likelihoods involves high-order multiple integrals) it may very time-consuming to obtain (by increasing the precision of the computation of the integral) a value $\tilde c$ below the value that we would have wished, for instance $10^{-2}$.

It may be interesting to stop at a value $\tilde c$ such that the approximation error is smaller than the statistical error but can not be considered as negligible, for instance $\tilde c= 0.5$. Let us denote $\tilde \theta$ a parameter value satisfying the convergence test with the stopping value $\tilde c$; $\tilde \theta$ is not far from $\hat \theta$ and in fact, considered as an estimator of $\theta$ it shares its consistency property. Moreover we can construct a conservative confidence region for $\theta$ based on $\tilde \theta$. An asymptotic $(1-\alpha$) confidence region for $\theta$ based on $\hat \theta$ can be defined as $R_{1-\alpha}=\{\theta: d(\theta,\hat \theta) \le c_{\alpha}\}$ where $d(\theta,\hat \theta)=U(\theta)^TG^{-1}(\hat \theta)U(\theta)$ and $c_{\alpha}$ is the $(1-\alpha)$ quantile of the $\chi^2_m$ distribution (Knight, 2000). Consider now the confidence region centered on $\tilde \theta$: $\tilde R_{1-\alpha}=\{\theta: d(\theta,\tilde \theta) \le c_{\alpha}+d(\tilde \theta,\hat \theta)\}$. We have that $R_{1-\alpha} \subset \tilde R_{1-\alpha}$. The result can be proved by using the triangular inequality $d(\theta,\tilde \theta) \le d(\theta,\hat \theta)+ d(\hat \theta,\tilde \theta)$ to show that $d(\theta,\tilde \theta) \le c_{\alpha}+d(\tilde \theta,\hat \theta)\} \Longrightarrow d(\theta,\hat \theta) \le c_{\alpha}\}$. Thus $\tilde R_{1-\alpha}$ is a conservative $(1-\alpha$) confidence region for $\theta$. Confidence regions do not give directly usual confidence intervals for each component of $\theta$ (because confidence regions take multiplicity into account). However we go from $R_{1-\alpha}$ to $\tilde R_{1-\alpha}$ by inflating $G^{-1}$ by the factor $(c_{\alpha}+d(\tilde \theta,\hat \theta))/c_{\alpha}$; it is tempting to do the same for confidence intervals, at least in an exploratory phase. 

As an example a $0.95$ confidence ellipsoid for $m=10$ is obtained with $c_{\alpha}=18.3$. If we stop at $c_k=0.4$ we have $d(\tilde \theta,\hat \theta)\approx 4$; so the conservative confidence ellipsoid is $\tilde R_{0.95}=\{\theta: d(\theta,\tilde \theta) \le 22.31\}$ and the inflating factor is $1.22$; applying this to confidence intervals this would lead to intervals inflated by $\sqrt 1.22=1.1$.

\newpage
\begin{center}

Table 1 : Convergence statistics from the simulation study comparing the Marquardt and RVS algorithms for findin the the MLE in a mixed model from left-censored longitudinal data.

\vspace{5mm}

\begin{tabular}{l c c}
\hline
  & Marquardt & RVS \\
\hline
Model AR & & \\
\hspace{2mm} threshold =$1.0$ & & \\
 \# convergence reached & 97 & 97 \\
 Mean convergence time & 47.7 & 21.0 \\
 Mean iterations number & 6.2 & 9.6 \\
\hspace{2mm} threshold =$2.0$ & & \\
 \# convergence reached & 83 & 97\\
 Mean convergence time & 206.3& 71.8 \\
 Mean iterations number & 9.6 & 11.5\\
 \hline
Model RE & & \\
\hspace{2mm} threshold =$1.0$ & & \\
 \# convergence reached & 95 & 98 \\
 Mean convergence time & 19.6 & 9.1\\
 Mean iterations number & 7.3 & 11.4\\
\hspace{2mm} threshold =$2.0$ & & \\
 \# convergence reached & 86 & 97\\
 Mean convergence time & 128.1 & 43.8\\
 Mean iterations number & 11.3 & 13.1\\
\hline
\end{tabular}
\end{center}

\newpage
\begin{center}

Table 2 : Asymptotic and sample variances of the estimates and coverage rates of the 95\% confidence intervals of the maximum likelihood estimators obtained with Marquardt and RVS algorithms when the censoring threshold was $2.0$.

\vspace{5mm}

\begin{tabular}{c c c c c c c }
\hline
  & \multicolumn{3}{c}{Marquardt} &  \multicolumn{3}{c}{RVS} \\
\hline
 & Asymptotic & Sample & Coverage & Asymptotic & Sample & Coverage \\
parameter & variance &  variance & rate  & variance &  variance & rate \\
Model AR  &  & & & & &  \\
$\beta_0$ &0.036 &0.040 &91.6 &0.038 & 0.040& 91.8\\
$\beta_1$ & 0.0021& 0.0027& 89.2& 0.0024&0.0024 & 91.8\\
$\beta_2$ & 0.072& 0.096& 91.6& 0.079& 0.091& 93.8\\
$\beta_3$ & 0.0048& 0.0056& 92.8& 0.0054& 0.0050& 94.9\\
$\sigma^2_w$ & 0.029&0.036 &88.0 &0.036 &0.034 & 92.8 \\
$\delta$ & 0.0027& 0.0027&90.4 & 0.0031&0.0025 & 94.9\\
$\sigma^2_e$ &0.0094 & 0.0095& 95.2& 0.010&0.009 &96.9 \\
\hline
Model RE  & & & & & & \\
$\beta_0$ & 0.018 &0.016 &94.2 &0.020 & 0.018& 94.9\\
$\beta_1$ & 0.0036& 0.0035& 95.4& 0.0041& 0.0034& 97.9\\
$\beta_2$ & 0.037 &0.037 & 97.7& 0.042&0.039 & 99.0\\
$\beta_3$ & 0.0077&0.0068 &96.5 &0.0088 &0.0067 & 99.0\\
$\sigma^2_0$ & 0.016&0.017 & 90.7& 0.021& 0.018& 94.9 \\
$\sigma_{01}$& 0.0026& 0.0026& 91.9& 0.0034&0.0030 &93.8 \\
$\sigma^2_{1}$ & 0.0008& 0.0008& 91.9& 0.0011& 0.0008&94.9 \\
$\sigma^2_e$ &0.0065 & 0.0070& 95.3& 0.0074&0.0064 & 95.9\\
\hline
\end{tabular}
\end{center}

\end{document}